\documentclass[10pt,reqno]{amsart}

\usepackage{amsmath,amsthm,amssymb,tikz-cd,comment,hyperref,thmtools,cleveref,nicematrix,enumerate,subcaption,float}

\theoremstyle{plain}
\newtheorem{thm}{Theorem}[section]
\newtheorem{lem}[thm]{Lemma}

\newtheorem{cor}[thm]{Corollary}
\newtheorem{thmalphabetintro}{Theorem}

\newtheorem{thmalphabetmaintext}{Theorem}

\theoremstyle{definition}

\newtheorem{exmp}[thm]{Example}
\newtheorem{rem}[thm]{Remark}

\newtheorem{notation}[thm]{Notation}

\newtheorem{remark and notation}[thm]{Remark and Notation}

\def \p {\mathbb{P}}

\def \N {\textup{\bf{N}}}
\def \H {\textup{H}}
\def \sing {\textup{Sing}}

\def \sec {\operatorname{Sec}}

\def \depth {\operatorname{depth}}
\def \tor {\operatorname{Tor}}
\def \reg {\operatorname{reg}}

\def\p{{\mathbb P}}

\title[Projective varieties beyond del Pezzo varieties]{Characterization of projective varieties beyond varieties of minimal degree and del Pezzo varieties}
\author{Jong In Han}
\address{Jong In Han, Department of Mathematical Sciences, Korea Advanced Institute of Science and Technology (KAIST), 291, Daehak-ro, Yuseong-gu, Daejeon, Republic of Korea}
\email{jihan09@kaist.ac.kr}

\author{Sijong Kwak}
\address{Sijong Kwak, Department of Mathematical Sciences, Korea Advanced Institute of Science and Technology (KAIST), 291, Daehak-ro, Yuseong-gu, Daejeon, Republic of Korea}
\email{sjkwak@kaist.ac.kr}

\author{Euisung Park}
\address{Euisung Park, Department of Mathematics, Korea University, 145, Anam-ro, Seongbuk-gu, Seoul, Republic of Korea}
\email{euisungpark@korea.ac.kr}

\date{}
\thanks{J. I. Han was partially supported by Basic Science Research Program through the National Research Foundation of Korea (NRF) funded by the Ministry of Education (2019R1A6A1A10073887). S. Kwak was partially supported by the National Research Foundation of Korea (NRF) grant funded by the Korea government (MSIT) (No. 2021R1A2C1013851). E. Park was partially supported by the National Research Foundation of Korea (NRF) grant funded by the Korea government (MSIT) (No. 2022R1A2C1002784).}

\begin{document}
\begin{abstract}
	Varieties of minimal degree and del Pezzo varieties are basic objects in projective algebraic geometry. Those varieties have been characterized and classified for a long time in many aspects. Motivated by the question ``which varieties are the most basic and simplest except the above two kinds of varieties in view of geometry and syzygies?", we give an upper bound of the graded Betti numbers in the quadratic strand and characterize the extremal cases.
	
	The extremal varieties of dimension $n$, codimension $e$, and degree $d$ are exactly characterized by the following two types:
	\begin{enumerate}
		\item[(i)] Varieties with $d = e+2$, $\depth X =n$, and Green-Lazarsfeld index $a(X)=0$,
		\item[(ii)] Arithmetically Cohen-Macaulay varieties with $d = e+3$.
	\end{enumerate}
	This is a generalization of G. Castelnuovo, G. Fano, and E. Park’s results on the number of quadrics and an extension of the characterizations of varieties of minimal degree and del Pezzo varieties in view of linear syzygies of quadrics due to K. Han and S. Kwak (\cite {{zbMATH02692308},{zbMATH02683130},{MR3334084},{MR3302629}}).
	
	In addition, we show that every variety $X$ that belongs to (i) or (ii) is always contained in a unique rational normal scroll $Y$ as a divisor. Also, we describe the divisor class of $X$ in $Y$.
\end{abstract}	

\maketitle
\tableofcontents
\section{Introduction}
Let $X^n\subseteq\p^{n+e}$ be a nondegenerate projective variety of dimension $n$, codimension $e$, and degree $d$ over an algebraically closed field $k$ of characteristic 0. We denote the $(p,q)$-th Betti number of $X$ as $\beta_{p,q}(X):=\dim\tor_p(S_X,k)_{p+q}$ where $S_X$ is the homogeneous coordinate ring of $X$. We also denote the saturated ideal of $X$ as $I_X$.
We say $X$ satisfies property $\N_{2,p}$ if $\beta_{i,j}(X)=0$ for all $i\leq p$ and $j\geq 2$, i.e. $S_X$ has only linear syzygies up to $p$-th step in the minimal free resolution. Let $a(X)$ be the largest integer $p$ (or $\infty$) such that $X$ satisfies property $\N_{2,p}$ which is called \textit{the Green-Lazarsfeld index} (\textit{c.f.} \cite{MR2854112}).
The depth of $X$ is defined as the depth of its homogeneous coordinate ring $S_X$ and is denoted as $\depth X$.
It is well known that $1\leq \depth X\leq n+1$, and we say X is arithmetically Cohen-Macaulay (ACM) if $\depth X=n+1$, i.e. $S_X$ is Cohen-Macaulay.

It holds that $d\geq e+1$, and $X$ is called a variety of minimal degree if $d=e+1$. Varieties of minimal degree are always normal and ACM. Also $X$ is called a del Pezzo variety if it is an ACM variety of degree $d=e+2$.
Varieties of minimal degree and del Pezzo varieties are the most basic objects in projective algebraic geometry, and they have been characterized and classified for a long time.
Classically, the upper bounds of the number of quadrics $\beta_{1,1}(X)$ containing $X$ and the extremal cases of the bounds were studied by Castelnuovo and Fano:
\\
\\
\noindent
\textbf{(a)}
\textit{
	Let $X^n\subseteq\p^{n+e}$ be a nondegenerate projective variety of codimension $e$. Then we have
	\[
		\beta_{1,1}(X)\leq \dbinom{e+1}{2},
	\]
	and the equality holds if any only if $X$ is a variety of minimal degree \textup{(\cite{zbMATH02692308})}.
	Also, unless $X$ is a variety of minimal degree, we have
	\[
		\beta_{1,1}(X)\leq \dbinom{e+1}{2} -1,
	\]
	and the equality holds if any only if $X$ is a del Pezzo variety \textup{(\cite{zbMATH02683130})}.
}
\\
\\
\indent
F. Zak reproved the results of Castelnuovo and Fano in \cite{MR1678545}. There is also another direction on the bounds of the number of equations: J. Harris and S. L'vovsky gave the bounds of the number of $m$-forms $h^0(\p^{e+1},\mathcal{I}_X(m))$ when $X$ is a curve (\cite{MR616900},\cite{MR1418349}).

We say a Betti number $\beta_{p,q}$ is in the quadratic strand if $q=1$. The Betti numbers in the quadratic strand have been interesting projective invariants after M. Green's fundamental works including the $K_{p,1}$-theorem (\textit{c.f.} Theorem \ref{K_p,1}). As a generalization of results of Castelnuovo and Fano, the upper bounds of Betti numbers in the quadratic strand and the extremal cases have been known as follows:
\\
\\
\noindent
\textbf{(b)} (\cite[Theorem 1.2, Theorem 1.3]{MR3302629})
\textit{
	Let $X^n\subseteq\p^{n+e}$ be a nondegenerate projective variety of codimension $e$. Then we have
	\[
		\beta_{p,1}(X)\leq
		\begin{cases}
			p\dbinom{e+1}{p+1} &\text{if } 1\leq p\leq e\\
			0 &\text{if } p\geq e+1,\\
		\end{cases}
	\]
	and the following are equivalent:
	\begin{enumerate}
		\item[$(i)$] $X$ is a variety of minimal degree, i.e. $d=e+1$.
		\item[$(ii)$] $\beta_{p,1}(X)=p\dbinom{e+1}{p+1}$ for some $1\leq p\leq e$.\\
		\item[$(iii)$] $\beta_{p,1}(X)=p\dbinom{e+1}{p+1}$ for all $1\leq p\leq e$.
		\item[$(iv)$] $X$ satisfies property $\N_{2,e}$.
		\item[$(v)$] $X$ has the following Betti table
		\[\begin{array}{c|ccccc}
			&0 &1 &2 &\cdots &e \\ \hline
			0&1 & & & & \\
			1& &\beta_{1,1} &\beta_{2,1} &\cdots &\beta_{e,1}
		\end{array}\]
		where $\beta_{p,1}=p\dbinom{e+1}{p+1}$ for all $1\leq p\leq e$.
	\end{enumerate}
}

By the characterizations above, the varieties of minimal degree are considered as the simplest projective varieties in syzygetic viewpoints. As a next step, for the varieties that are not of minimal degree, the following upper bounds of Betti numbers in the quadratic strand and the extremal cases are known:
\\
\\
\noindent
\textbf{(c)} (\cite[Theorem 1.4, Theorem 4.3]{MR3302629})
\textit{
	Let $X^n\subseteq\p^{n+e}$ be a nondegenerate projective variety of codimension $e\geq 2$. Unless $X$ is a variety of minimal degree, we have
	\[
		\beta_{p,1}(X)\leq
		\begin{cases}
			p\dbinom{e+1}{p+1}-\dbinom{e}{p-1} &\text{if } 1\leq p\leq e-1\\
			0 &\text{if } p\geq e,\\
		\end{cases}
	\]
	and the following are equivalent:
	\begin{enumerate}
		\item[$(i)$] $X$ is a del Pezzo variety, i.e. $X$ is an ACM variety of degree $d=e+2$.
		\item[$(ii)$] $\beta_{p,1}(X)=p\dbinom{e+1}{p+1}-\dbinom{e}{p-1}$ for some $1\leq p\leq e-1$.\\
		\item[$(iii)$] $\beta_{p,1}(X)=p\dbinom{e+1}{p+1}-\dbinom{e}{p-1}$ for all $1\leq p\leq e-1$.
		\item[$(iv)$] $X$ satisfies property $\N_{2,e-1}$, but not property $\N_{2,e}$.
		\item[$(v)$] $X$ has the following Betti table
		\[\begin{array}{c|ccccc}
			&0 &1 &\cdots &e-1 &e \\ \hline
			0&1 & & & & \\
			1& &\beta_{1,1} &\cdots &\beta_{e-1,1} &\\
			2& & & & &1
		\end{array}\]
		where $\beta_{p,1}=p\dbinom{e+1}{p+1}-\dbinom{e}{p-1}$ for all $1\leq p\leq e-1$.
	\end{enumerate}
}
~\\
\indent
From this point of view, del Pezzo varieties are considered as the next simplest projective varieties.
Then it is natural to ask what are the next basic objects in this direction. To answer the question, we first induce the next upper bounds of Betti numbers:

\begin{thmalphabetintro}\label{intro_upper_bound_not_vmd_dp}
	Let $X^n\subseteq\p^{n+e}$ be a nondegenerate projective variety of codimension $e\geq 2$. Unless $X$ is either a variety of minimal degree or a del Pezzo variety, we have
	\[
	\beta_{p,1}(X)\leq
	\begin{cases}
		p\dbinom{e+1}{p+1}-2\dbinom{e}{p-1} &\text{if } 1\leq p\leq e-1\\
		0 &\text{if } p\geq e.\\
	\end{cases}
	\]
\end{thmalphabetintro}

Then we characterize the extremal cases as follows:

\begin{thmalphabetintro}\label{intro_extremal_case}
	Let $X^n\subseteq\p^{n+e}$ be a nondegenerate projective variety of dimension $n$, codimension $e\geq 3$, degree $d$, and Green-Lazarsfeld index $a(X)$.
	Then the following are equivalent:
	{\small
	\begin{enumerate}
		\item[$(i)$] $X$ is either a variety of almost minimal degree (\textit{i.e.} $d=e+2$) with $\depth X=n$ and $a(X)=0$, or an ACM variety of degree $d=e+3$.
		\item[$(ii)$] $\beta_{p,1}(X)=p\dbinom{e+1}{p+1}-2\dbinom{e}{p-1}$ for an integer $\begin{cases}
			p=1 &\text{if }e=3\\
			2\leq p\leq e-2 &\text{if }e\geq 4.
		\end{cases}$
		\item[$(iii)$] $\beta_{p,1}(X)=p\dbinom{e+1}{p+1}-2\dbinom{e}{p-1}$ for all $1\leq p\leq e-1$.
	\end{enumerate}
	}
\end{thmalphabetintro}

\Cref{intro_extremal_case} shows that there are two types of extremal varieties with respect to \Cref{intro_upper_bound_not_vmd_dp}. The second type, i.e. the ACM varieties of degree $d=e+3$, would be the naturally expected one according to the characterization result in \cite{MR3302629}.
One can get a parallel statement to \textbf{(a)} and \textbf{(b)} for the ACM varieties of degree $d=e+3$ (\Cref{parallel_statement_d=e+3}).
On the other hand, the first case is unexpected. The existence of this case, or the exclusion of the other cases, can be shown using the inner projection method (\Cref{exclude_other_depth_n_VAMD}).
Also, it can be explained geometrically why these two types of extremal varieties in \Cref{intro_extremal_case} have the same Betti numbers in the quadratic strand.
Roughly speaking, the union of a variety of the first type and an $n$-dimensional linear space is a reducible ACM algebraic set of degree $e+3$, which has the same Betti table as the varieties of the second type. For details, see \Cref{same_quadratic_strand}.

The two extremal varieties in \Cref{intro_extremal_case} are always divisors in some rational normal scrolls (\textit{c.f.} Section 4).
Along this line, we investigate the numerical type of such rational normal scrolls and also the divisor classes of the extremal varieties.

\begin{thmalphabetintro}\label{intro_VAMD_divisor}
Let $X^n \subseteq \p^{n+e}$ be a nondegenerate projective variety of dimension $n$, codimension $e \geq 3$, and degree $d$. Then the following two statements are equivalent:
\begin{enumerate}
	\item[$(i)$] X is a variety of almost minimal degree with $\depth X=n$ and $a(X)=0$.
	\item[$(ii)$] $X$ is contained in an (unique) $(n+1)$-fold rational normal scroll $Y = S(0,\ldots , 0,1,e-1)$ as a divisor linearly equivalent to $H+2F$ where $H$ and $F$ are respectively the hyperplane section and a ruling of $Y$.
\end{enumerate}
Furthermore, in the case above, if $X$ is not a cone then $n \leq 4$.
\end{thmalphabetintro}

\begin{thmalphabetintro}\label{intro_ACM_e+3_divisor}
	Let $X^n \subseteq \p^{n+e}$ be a nondegenerate projective variety of dimension $n$, codimension $e \geq 3$, and degree $d$. Then the following two statements are equivalent:
	\begin{enumerate}
		\item[$(i)$] $X$ is an ACM variety of degree $d=e+3$.
		\item[$(ii)$] $X$ is contained in an (unique) $(n+1)$-fold rational normal scroll $Y = S(a_0 , a_1 , \ldots ,a_n )$ for some $0 \leq a_0 \leq a_1 \leq \ldots \leq a_n$ such that either
		\begin{enumerate}
			\item[$(\alpha)$] $a_{n-1} > 0$ and the divisor class of $X$ is $2H+(3-e)F$  where $H$ and $F$ are respectively the hyperplane section and a ruling of $Y$, or
			\item[$(\beta)$] $a_{n-1} = 0$ and the divisor class of $X$ is $(e+3)R$  where $R \cong \p^n$ is the effective generator of the divisor class group of $Y$.
		\end{enumerate}
	\end{enumerate}
\end{thmalphabetintro}

The classification of ACM varieties of degree $d=e+3$ has not been completed yet, but for the smooth case, P. Ionescu completely classified them over $\mathbb{C}$ (\cite[Theorem 3.4]{MR749942}, also \textit{c.f.} \cite[(10.7)]{MR1162108}). Interestingly, his result shows smooth nondegenerate ACM projective varieties of degree $d=e+3$ must have dimension $\leq 4$. Furthermore, when $n=4$, the variety must be the Segre product of $\p^1$ and a smooth hyperquadric $Q^3\subseteq\p^4$ while there are five cases for $n=3$.
Before Ionescu classified such varieties of dimension $\geq 3$, Castelnuovo also studied the case of surfaces in \cite{zbMATH02690412} (\textit{c.f.} \cite[Proposition 3.1]{MR749942}, also \cite[(10.7.2)]{MR1162108}).
~\\
\\
\indent
We review the well-known facts used in the proofs of our main theorems in Section 2. In Section 3, we prove the main theorems and provide some examples supporting our results. In Section 4, we describe the unique rational normal scroll containing an extremal variety as a divisor so as to prove \Cref{intro_VAMD_divisor} and \Cref{intro_ACM_e+3_divisor}. We used Macaulay2 (\cite{M2})
to check the Betti numbers of various projective varieties throughout this paper. It was
also useful to estimate the upper bounds.

\section{Preliminaries}

\noindent
\textit{2.1.}\textbf{ Definitions and fundamental results.}
\begin{notation}
	In addition to the notations from the introduction, we introduce a few more notations and conventions we use.
	\begin{itemize}
		\item When we say a variety, we mean an integral projective variety over an algebraically closed field of characteristic zero.
		\item A variety of codimension $e$ and degree $d$ is called a \textit{variety of minimal degree} (VMD) when $d=e+1$ and called a \textit{variety of almost minimal degree} when $d=e+2$.
		\item The singular locus of a variety $X$ is denoted as $\sing(X)$.
		\item The join of two varieties $X$ and $Y$ in $\p^r$ is denoted as $J(X,Y)$.
		\item The vertex set of a variety $X\subseteq\p^{r}$ is the set $\{p\in X\mid J(p,X)=X\}$.
		\item The sectional genus $\pi(X)$ of a variety $X\subseteq\p^{r}$ is defined as the arithmetic genus of the general curve section of $X$.
		\item The depth of $X$ is defined as $\depth S_X$ which is equal to
		\[
			\depth X=\min\{i\in \mathbb{Z}_{\geq 1}\mid \H^i(\p^{n+e},\mathcal{I}_X(m))\neq 0\text{ for some }m\in\mathbb{Z}\}.
		\]
		\item We denote the rational normal scroll of type $(b_1,\cdots,b_m)$ as $S(b_1,\cdots,b_m)$. We always assume the indices are ordered such that $0\leq b_1\leq \cdots \leq b_m$.
		\item For a projective variety $X^n\subseteq\p^{n+e}$ of dimension $n$ and codimension $e$, we define the regularity $\reg X$ as the regularity of its ideal sheaf $\mathcal{I}_X$ which is defined as
		\[
			\reg \mathcal{I}_X:=\min\{d\in\mathbb{Z}\mid \H^i(\p^{n+e},\mathcal{I}_X(d-i))=0\text{ for all }i\geq 1\}.
		\]
	\end{itemize}
\end{notation}

\begin{rem}[Other characterizations of Betti numbers]
	Let $X^n\subseteq\p^{n+e}=\p(V^\vee)$ be a nondegenerate projective variety and $S$ be the polynomial ring $Sym(V)$. Recall that the Betti number $\beta_{p,q}(X)$ of $X$ is defined as $\dim\tor_p(S_X,k)_{p+q}$ where $S_X$ is the homogeneous coordinate ring of $X$. The Betti numbers can be also characterized as follows:
	\begin{itemize}
		\item For a minimal free resolution,
		\[
		0\to F_m\to \cdots \to F_1 \to F_0 \to S_X
		\]
		of $S_X$, we have $F_p=\bigoplus\limits_{q}S(-p-q)^{\beta_{p,q}(X)}$.
		\item The Koszul cohomology group $K_{p,q}(M,V)$ is defined as the cohomology group at the middle term of the complex
		\[\wedge^{p+1}V\otimes M_{q-1}\to\wedge^{p}V\otimes M_{q}\to\wedge^{p-1}V\otimes M_{q+1}\]
		where the differential map $\wedge^{p}V\otimes M_{q}\to\wedge^{p-1}V\otimes M_{q+1}$ is given by $v_1\wedge\cdots\wedge v_p\otimes m\mapsto \sum\limits_{}^{}(-1)^i v_1\wedge\cdots\wedge \hat{v_i}\wedge\cdots\wedge v_p\otimes v_im$. Then $\beta_{p,q}(X)=\dim_k K_{p,q}(S_X,V)$.
	\end{itemize}
\end{rem}

The vanishing of Betti numbers and the distribution of nonzero Betti numbers of projective varieties have been considered important to study the geometry of projective varieties. In \cite{MR739785}, M. Green studied the Betti numbers of projective varieties thoroughly using Koszul cohomology and introduced the $K_{p,1}$-theorem which is one of the fundamental results in the subject.

\begin{thm}[{\cite[Theorem 3.c.1]{MR739785}, also \textit{c.f.} \cite[Theorem 3.5]{MR1291122}, $K_{p,1}$-theorem}]\label{K_p,1}
	Let $X^n\subseteq\p^{n+e}$ be a nondegenerate projective variety of dimension $n$, codimension $e\geq 2$, and degree $d$. Then
	\begin{enumerate}
		\item $\beta_{p,1}(X)=0$ for all $p>e$.
		\item $\beta_{e,1}(X)=0$ unless $X$ is a variety of minimal degree.
		\item $\beta_{e-1,1}(X)=0$ unless either $d\leq e+2$ or $X$ is contained in a variety of minimal degree of dimension $n+1$.
	\end{enumerate}
\end{thm}

\noindent
\textit{2.2.}\textbf{ Properties of inner projections.}
It is known that the general inner projections of projective varieties preserve some nice properties. Here, we introduce some theorems on inner projections we use.
The following inequality of Betti numbers plays a key role when we use the induction in the proofs of our main theorems.

\begin{thm}[{\cite[Theorem 3.1]{MR3302629}}]\label{betti_ineq}
	Let $X^n\subseteq\p^{n+e}$ be a nondegenerate projective variety of codimension $e$. Let $X_q\subseteq\p^{n+e-1}$ be the image of the inner projection of $X$ from a closed point $q\in X$. Then for $p\geq 1$, we have
	\[
	\beta_{p,1}(X)\leq \beta_{p,1}(X_q)+\beta_{p-1,1}(X_q)+\dbinom{e}{p},
	\]
	and equality holds if $q$ is a smooth point of $X$ and $1\leq p \leq a(X)$ where $a(X)$ is the largest integer (or $\infty$) such that $X$ satisfies property $\N_{2,a(X)}$. Furthermore, if the equality holds for some $p\geq 1$, then
	\[
	\beta_{1,1}(X)=\beta_{1,1}(X_q)+e.
	\]
\end{thm}

When $X$ is ideal-theoretically defined by quadrics, the general inner projection preserve the depth of $X$. Specifically,

\begin{thm}[{\cite[Theorem 4.1]{MR2946932}}]\label{depth_inner_projection}
	Let $X^n\subseteq\p^{n+e}$ be a nondegenerate reduced subscheme such that its saturated ideal $I_X$ is generated by quadrics. Let $X_q\subseteq\p^{n+e-1}$ be the image of the inner projection of $X$ from a smooth point $q\in X$. Then the arithmetic depths of $X\subseteq\p^{n+e}$ and $X_q\subseteq\p^{n+e-1}$ are same.
\end{thm}

As a result, we get the following lemma which is useful when we use the induction.

\begin{lem}\label{inner_projection_dp}
	Let $X^n\subseteq\p^{n+e}$ be a nondegenerate projective variety, and $X_q\subseteq\p^{n+e-1}$ be the image of the inner projection of $X$ from the general point $q\in X$. Then $X$ is a del Pezzo variety if and only if $X_q$ is a del Pezzo variety.
\end{lem}
\begin{proof}
	First note that $\deg X_q=\deg X -1$ as the general point $q\in X$ is a smooth point.
	If $X$ is del Pezzo, then $\depth X=\depth X_q$ by \Cref{depth_inner_projection}. Hence $X_q$ is also del Pezzo. Conversely, suppose $X_q$ is a del Pezzo variety while $X$ is not. Let $C\subseteq\p^{e+1}$ be the general curve section of $X$ containing $q$ and $C_q\subseteq\p^{e}$ be the inner projection of $C$ at $q$. 
	Then $C_q$ is a general curve section of $X_q$.
	Also, $C_q$ is a del Pezzo curve while $C$ is not. In this case, it is well known that $C$ is a smooth rational curve of degree $e+2$ and the inner projection of $C$ at $q$ is an embedding since there is no trisecant line through $q$. In particular, $C$ and $C_q$ should have the same genus. This is a contradiction.
\end{proof}

\noindent
\textit{2.3.}\textbf{ General hyperplane sections and syzygies of finite points.}
Before the inner projections were studied, the classical way to reduct the projective varieties is to take the general hyperplane sections.
Recall that for a nondegenerate projective variety $X^n\subseteq\p^{n+e}$ of depth $t\geq 2$, the general hyperplane section of $X$ has the same Betti numbers with $X$.
The following theorem, which is known as the Lefschetz theorem, enables us to take a hyperplane section of $X$ even if $\depth X=1$ when we consider the upper bound of Betti numbers.

\begin{thm}[{\textit{c.f.} \cite[Theorem 2.20]{MR2573635}, Lefschetz theorem}]\label{Lefschetz}
	Let $X^n\subseteq\p^{n+e}$ be a nondegenerate projective variety of dimension $n$ and $Y^{n-1}\subseteq\p^{n+e-1}$ be a hyperplane section of $X$. Then $\beta_{p,1}(X)\leq \beta_{p,1}(Y)$ where $\beta_{p,1}(Y)$ be the $(p,1)$-th Betti number of $Y$ as a projective variety in $\p^{n+e-1}$.
\end{thm}

With the Lefschetz theorem, the syzygies of finite points provide useful information when we get the upper bound of Betti numbers and characterize the extremal cases. The following theorems are the results on the syzygies of finite points we use.

\begin{thm}[{\cite[Theorem 3.c.6]{MR739785}, Strong Castelnuovo lemma}]\label{strong_Castelnuovo}
	Let $X\subseteq\p^r$ be a set of $d$ points in general position. Then $X$ lies on a rational normal curve if and only if $\beta_{r-1,1}(X)\neq 0$.
\end{thm}

\begin{thm}[{\cite[Lemma 2.1]{MR1291122}}]\label{finite_pts_on_rnc}
	Let $X\subseteq\p^r$ be a finite scheme of degree $d$ on a rational normal curve $C\subseteq\p^r$. Let $p$ be an integer such that $0\leq p\leq n$. Then $\beta_{p,1}(X)=\beta_{p,1}(C)$ for all $p\geq 2r+2-d$.
\end{thm}

\begin{thm}[{\cite[Theorem 1]{MR959214}}]\label{finite_pts_Np}
	Let $X\subseteq\p^r$ be a set of $2r+1-p$ points in general position. Then $X$ satisfies property $\N_{2,p}$.
\end{thm}

\noindent
\textit{2.4.}\textbf{ Varieties of almost minimal degree.}
Now we briefly review the classification of the Betti tables of varieties of almost minimal degree that are not ACM.

Let $X^n\subseteq\p^{n+e}$ be a variety of almost minimal degree with $\depth X=t$ where $1\leq t\leq n$.
Let $\bar{X}\subseteq\p^{n+e-t+1}$ be the general linear section of dimension $n-t+1$ of $X$.
Then $X$ and $\bar{X}$ have the same Betti table.
In particular, when $e\geq 3$, $\bar{X}$ is contained in a rational normal scroll $S(b_1,\cdots,b_{n-t+2})$ as a divisor for some integers $1\leq b_1\leq \cdots\leq b_{n-t+2}$ such that $b_1+\cdots+b_{n-t+2}=e$ (\cite[Theorem 1.1]{MR2331758}).
Furthermore, $\bar{X}$ is linearly equivalent to $H+2F$ where $H$ is the hyperplane section of $S(b_1,\cdots,b_{n-t+2})$ and $F$ is the fiber of the map $S(b_1,\cdots,b_{n-t+2})\to \p^1$.

We denote the Betti table of the divisor $H+2F$ of $S(b_1,\cdots,b_{n-t+2})$ as $T(b_1,\cdots,b_{n-t+2})$ following the notation of \cite{MR2850075}.
Note that the Betti tables of the divisors in such class are all the same so that $\bar{X}$ has the Betti table $T(b_1,\cdots,b_{n-t+2})$.
For the case $e=2$, there are only two cases: $\bar{X}$ is a divisor of $S(1,1)$ that is linearly equivalent to $H+2F$ whose Betti table is denoted as $T(1,1)$, or $\bar{X}\subseteq\p^4$ is a projection of the Veronese surface $\nu_2(\p^2)\subseteq\p^5$ from a point outside of the secant variety of $\nu_2(\p^2)$, and we denote the Betti table as $T(V)$.
Consequently, we can summerize all possible Betti tables as follows:

\begin{thm}[{\cite[Theorem 3.1]{MR2850075}}]\label{classification_Betti_VAMD}
	Let $X^n\subseteq\p^{n+e}$ be a nondegenerate projective variety of almost minimal degree of dimension $n$, codimension $e\geq 2$, and $\depth X=t$ $(1\leq t\leq n)$.
	\begin{enumerate}
		\item When $e=2$, the Betti table of $X$ is equal to either $T(V)$ or $T(1,1)$.
		\item When $e\geq 3$, the Betti table of $X$ is equal to $T(b_1,\cdots,b_{n-t+2})$ for some integers $1\leq b_1\leq \cdots\leq b_{n-t+2}$ such that $b_1+\cdots+b_{n-t+2}=e$.
	\end{enumerate}
\end{thm}

Although not all of the Betti numbers of varieties of almost minimal degree can be expressed explicitly, some Betti numbers including $\beta_{1,1}$, and $\beta_{e-1,1}$ can be:

\begin{thm}[{\cite{MR1228155}, \cite{MR1615938}, \cite{MR2274517}}]\label{betti_almost_minimal}
	Let $X^n\subseteq\p^{n+e}$ be a nondegenerate projective variety of almost minimal degree with codimension $e\geq 2$ and $\depth X=t$. If $X$ is not del Pezzo, then \begin{align*}
		\beta_{1,1}(X) &= \dbinom{e+1}{2}+t-n-2,\\
		\beta_{e-1,1}(X) &= e-1,\\
		\beta_{p,1}(X)-\beta_{p-1,2}(X) &= (e+1)\dbinom{e}{p}-\dbinom{e}{p+1}-\dbinom{n+e-t+2}{p}
	\end{align*}
	for all $2\leq p\leq e$.
\end{thm}

On the other hand, generalizing Castelnuovo and Fano's results, we have the following:

\begin{thm}[{\cite[Corollary 1.4]{MR3334084}, Characterization of $\beta_{1,1}(X)=\binom{e+1}{2}-2$}]\label{characterization_quadrics=max-2}
	Let $X^n\subseteq\p^{n+e}$ be a nondegenerate projective vairety of dimension $n$, codimension $e\geq 3$, and degree $d$. Then $$\beta_{1,1}(X)=\dbinom{e+1}{2}-2$$ if and only if either $d=e+2$ and $\depth X=n$, or $d=e+3$ and $\depth X=n+1$.
\end{thm}

\section{Proofs of the main theorems}
In this section, we give the sharp upper bounds of Betti numbers in the quadratic strand when a given variety is neither a variety of minimal degree nor a del Pezzo variety, and give characterizations of the extremal cases.
\newline
\newline
\noindent
\textit{3.1.} \textbf{Sharp upper bounds of Betti numbers in the quadratic strand.}

\begin{thmalphabetmaintext}[Sharp upper bounds of Betti numbers in the quadratic strand]\label{main_upper_bound_not_vmd_dp}
	Let $X^n\subseteq\p^{n+e}$ be a nondegenerate projective variety of codimension $e$. Unless $X$ is either a variety of minimal degree or a del Pezzo variety, we have
	\[
	\beta_{p,1}(X)\leq
	\begin{cases}
		p\dbinom{e+1}{p+1}-2\dbinom{e}{p-1} &\text{if } 1\leq p\leq e-1\\
		0 &\text{if } p\geq e.\\
	\end{cases}
	\]
\end{thmalphabetmaintext}
\begin{proof}
	We use the induction on $e$. For the case $e=2$, it is trivial to show $\beta_{1,1}(X)\leq 1$ since $X$ is neither of minimal degree nor del Pezzo.
	Now let $X$ be a nondegenerate projective variety of codimension $e\geq 3$, and suppose the inequality holds for varieties of codimension $\leq e-1$. It is immediate $$\beta_{1,1}(X)\leq \dbinom{e+1}{2}-2$$ by the result of Fano. Denote the image of the inner projection of $X$ from the general point $q\in X$ as $X_q\subseteq\p^{n+e-1}$.
	Note that if $X$ is neither a variety of minimal degree nor a del Pezzo variety, then so is $X_q$ by \Cref{inner_projection_dp}.
	Thus for any $2\leq p\leq e-2$, we have
	\begin{align*}
		\beta_{p,1}(X) &\leq \beta_{p,1}(X_q)+\beta_{p-1,1}(X_q)+\dbinom{e}{p}\\
		&\leq p\dbinom{e}{p+1}-2\dbinom{e-1}{p-1}+(p-1)\dbinom{e}{p}-2\dbinom{e-1}{p-2}+\dbinom{e}{p}\\
		&=p\dbinom{e+1}{p+1}-2\dbinom{e}{p-1}
	\end{align*}
	by \Cref{betti_ineq}.
	
	For $p=e-1$, we can not use the induction argument above.
	Suppose $\beta_{e-1,1}(X)\neq 0$. We divide into two cases to show $\beta_{e-1,1}(X)\leq e-1$. For the first case, suppose $d\geq e+3$. Note that the general $e$-plane section $Y$ of $X$ consists of $d$ points in general position. Let $\beta_{e-1,1}(Y)$ be the Betti number of $Y$ where $Y$ is considered as an algebraic set in $\p^e$. As $\beta_{e-1,1}(Y)\geq\beta_{e-1,1}(X)>0$ by the Lefschetz theorem (\Cref{Lefschetz}), $Y$ lies on a rational normal curve $C\subseteq\p^e$ by the strong Castelnuovo lemma (\Cref{strong_Castelnuovo}). Then by \Cref{finite_pts_on_rnc}, we get $\beta_{e-1,1}(Y)=\beta_{e-1,1}(C)=e-1$. Consequently, $\beta_{e-1,1}(X)\leq e-1$.
	
	For the second case, suppose $d=e+2$ and $X$ is not del Pezzo. Then it is done as $\beta_{e-1,1}(X)=e-1$ by \Cref{betti_almost_minimal}.
\end{proof}

From \Cref{main_upper_bound_not_vmd_dp}, we get $\beta_{e-1,1}(X)\leq e-1$ if $X$ is neither a variety of minimal degree nor a del Pezzo variety. However, $\beta_{e-1,1}(X)$ is a very rigid number in the sense that there are only two possible values of it actually:

\begin{cor}\label{beta_e-1,1}
	Let $X^n\subseteq\p^{n+e}$ be a nondegenerate projective variety of codimension $e$ and degree $d$. Unless $X$ is a variety of minimal degree or a del Pezzo variety, we have $\beta_{e-1,1}(X)=0$ or $e-1$.
\end{cor}
\begin{proof}
	If $X$ is of almost minimal degree that is not del Pezzo, then the result follows from \Cref{betti_almost_minimal}. Hence we assume $d\geq e+3$. Then $X$ is contained in a variety of minimal degree $Y$ as a divisor by the $K_{p,1}$-theorem so that $\beta_{e-1,1}(X)\geq\beta_{e-1,1}(Y)=e-1$. Therefore the result follows by \Cref{main_upper_bound_not_vmd_dp}.
\end{proof}

\begin{rem}
	For a nondegenerate projective variety $X^n\subseteq\p^{n+e}$ of codimension $e$, we can summerize the values of two highest Betti numbers in the quadratic strand as follows:
	\begin{align*}
		\beta_{e,1}(X) &=
		\begin{cases}
			e &\text{if $X$ is a VMD, i.e. a variety of minimal degree}\\
			0 &\text{otherwise}
		\end{cases}
		\\		
		\beta_{e-1,1}(X) &=
		\begin{cases}
			e^2-1 &\text{if $X$ is a VMD}\\
			\binom{e+1}{2}-1 &\text{if $X$ is a del Pezzo variety}\\
			e-1 &\text{else if $X$ is a divisor of a VMD}\\
			0 &\text{otherwise}
		\end{cases}
	\end{align*}
\end{rem}
~\\
\\
\noindent
\textit{3.2.} \textbf{Characterizations of extremal cases.}
Next, we characterize the extremal cases of the upper bounds. The extremal cases have $$\beta_{1,1}(X)=\dbinom{e+1}{2}-2$$ so that $X$ is an ACM variety of degree $d=e+3$ or a variety of almost minimal degree with $\depth X=n$ by \Cref{characterization_quadrics=max-2}.
When $X$ is a variety of almost minimal degree with codimension $e\geq 3$ and depth $n$, the Betti table of $X$ is $T(b_1,b_2)$ where $1\leq b_1\leq b_2$ and $b_1+b_2=e$.\footnote{See Section \textit{2.4} for the definition of $T(b_1,b_2)$.} The following lemma shows the Betti tables $T(b_1,b_2)$ can not be extremal if $b_1\neq 1$.

\begin{lem}\label{T(b_1,b_2)_bound}\label{exclude_other_depth_n_VAMD}
	Let $X^n\subseteq\p^{n+e}$ be a nondegenerate projective variety of almost minimal degree with $\depth X=n$ whose general $(e+1)$-plane section of $X$ is a smooth divisor of $S(b_1,e-b_1)$ for some $b_1\geq 2$. Then $$\beta_{p,1}(X)<p\dbinom{e+1}{p+1}-2\dbinom{e}{p-1}$$ for all $2\leq p\leq e-2$.
\end{lem}
\begin{proof}
	Suppose $b_1\geq 2$, and let $C\subseteq\p^{e+1}$ be the general $(e+1)$-plane section of $X$ that is a smooth divisor of $S(b_1,e-b_1)$ linearly equivalent to $H+2F$. Then $C$ and $X$ has the same Betti table $T(b_1,e-b_1)$ since $\depth X=n$. We proceed to show the inequality using the induction on $b_1$. When $b_1=2$, it is known that $$\beta_{p-1,2}(C)=\dbinom{e-1}{p-3}$$ for all $2\leq p\leq e-2$ by \cite[Theorem 1.1]{MR3463203}.\footnote{We use the convention $\binom{a}{b}=0$ for all $b<0$.}
	By \Cref{betti_almost_minimal}, we have
	\begin{align*}
		\beta_{p,1}(C) &= (e+1)\dbinom{e}{p}-\dbinom{e}{p+1}-\dbinom{e+2}{p}+\dbinom{e-1}{p-3}\\
		&< (e+1)\dbinom{e}{p}-\dbinom{e}{p+1}-\dbinom{e+2}{p}+\dbinom{e}{p-2}\\
		&= p\dbinom{e+1}{p+1}-2\dbinom{e}{p-1}
	\end{align*}
	for all $2\leq p\leq e-2$.
	Suppose $b_1\geq 3$, and the result holds until $b_1-1$. Take the general $(e+1)$-plane section $C\subseteq\p^{e+1}$ of $X$ which is a smooth curve contained in $S(b_1,e-b_1)$.
	Let $C_q\subseteq\p^e$ and $S_q\subseteq\p^e$ be the images of inner projections of $C$ and $S(b_1,e-b_1)$ from a point $q\in C\cap S(b_1)$, respectively. Then $S_q=S(b_1-1,e-b_1)$ and $C_q\subseteq S_q$. Since $C_q$ is a smooth curve of almost minimal degree with $\depth C_q=1$, by the induction hypothesis,
	\begin{align*}
		\beta_{p,1}(C) &\leq \beta_{p,1}(C_q)+\beta_{p-1,1}(C_q)+\dbinom{e}{p}\\
		&< p\dbinom{e}{p+1}-2\dbinom{e-1}{p-1}+(p-1)\dbinom{e}{p}-2\dbinom{e-1}{p-2}+\dbinom{e}{p}\\
		&= p\dbinom{e+1}{p+1}-2\dbinom{e}{p-1}.
	\end{align*}
\end{proof}

\begin{rem}
	Lemma \ref{T(b_1,b_2)_bound} improves the upper bounds of Betti numbers due to L. T. Hoa, M. Brodmann, and P. Schenzel (\cite[Theorem 2]{MR1228155}, \cite[Theorem 8.3.(a)]{MR2274517}). Specifically, they induced the bounds
	\[
		\beta_{p,1}(X)\leq (e+1)\dbinom{e}{p}-\dbinom{e}{p+1},\ 2\leq p\leq e-2
	\]
	for non-ACM varieties of almost minimal degree. The bounds in Lemma \ref{T(b_1,b_2)_bound} are smaller than these bounds by
	\[
		\left[(e+1)\dbinom{e}{p}-\dbinom{e}{p+1}\right]-\left[ p\dbinom{e+1}{p+1}-2\dbinom{e}{p-1}\right]=\dbinom{e}{p-1}+\dbinom{e+1}{p}.
	\]
\end{rem}

\begin{rem}
	The Betti table $T(b_1,b_2)$ satisfies $\N_{2,b_1-1}$ but not $\N_{2,b_1}$ for all $1\leq b_1\leq b_2$ by \cite[Corollary 4.2]{MR2299577}. In addition, the varieties of almost minimal degree of $\depth X=n$ and Green-Lazarsfeld index $a(X)=0$ (i.e. $b_1=1$) have the Betti table $T(1,e-1)$ which is equal to
	\[\begin{array}{c|cccccc}
		&0 &1 &\cdots &e-1 &e &e+1 \\ \hline
		0&1 & & & & & \\
		1& &\beta_{1,1} &\cdots &\beta_{e-1,1} & \\
		2& &\beta_{1,2} &\cdots &\beta_{e-1,2} &\beta_{e,2} &\beta_{e+1,2} 
	\end{array}\]
	\noindent
	where $\beta_{p,1}=p\dbinom{e+1}{p+1}-2\dbinom{e}{p-1}$ for all $1\leq p\leq e-1$ and
	\[\beta_{p,2}=
	\begin{cases}
		\dbinom{e}{p-1} &\text{if } 1\leq p\leq e-2\\[10pt]
		\dbinom{e+2}{p+1}-(e+1)\dbinom{e}{p+1} &\text{if } e-1\leq p\leq e+1
	\end{cases}\]
	by \cite[Theorem 1.1]{MR3463203}.
	More specifically, the general $(e+1)$-plane section of $X$ is a smooth divisor of the rational normal surface scroll $S(1,e-1)$ that is linearly equivalent to $H+2F$.
\end{rem}

\begin{thmalphabetmaintext}[Characterizations of extremal cases]\label{characterizations_extremal_cases}
	Let $X^n\subseteq\p^{n+e}$ be a nondegenerate projective variety of dimension $n$, codimension $e\geq 3$, degree $d$, and Green-Lazarsfeld index $a(X)$.
	Then the following are equivalent:
	{\small
		\begin{enumerate}
			\item[$(i)$] $X$ is either a variety of almost minimal degree (\textit{i.e.} $d=e+2$) with $\depth X=n$ and $a(X)=0$, or an ACM variety of degree $d=e+3$.
			\item[$(ii)$] $\beta_{p,1}(X)=p\dbinom{e+1}{p+1}-2\dbinom{e}{p-1}$ for an integer $\begin{cases}
				p=1 &\text{if }e=3\\
				2\leq p\leq e-2 &\text{if }e\geq 4.
			\end{cases}$
			\item[$(iii)$] $\beta_{p,1}(X)=p\dbinom{e+1}{p+1}-2\dbinom{e}{p-1}$ for all $1\leq p\leq e-1$.
		\end{enumerate}
	}
\end{thmalphabetmaintext}
\begin{proof}
	The Betti numbers of variety of almost minimal degree of $\depth X=n$ and $a(X)=0$ are completely determined by W. Lee and the third author in \cite[Theorem 1.1]{MR3463203}.
	Hence for $(i)\Rightarrow (iii)$, we only consider the case $X$ is an ACM variety of degree $d=e+3$.
	We use the induction on $e$ to show $(i)\Rightarrow (iii)$. Note that the general $e$-plane section $Y$ of $X$ is a set of $e+3$ points in general position in $\p^e$ whose Betti table is same as that of $X$. For the case $e=3$, we have $\beta_{1,1}(X)=4$ by \Cref{characterization_quadrics=max-2}, and $\beta_{2,1}(X)=2$ by \Cref{finite_pts_on_rnc} as the general $e$-plane section $Y$ is contained in a rational normal scroll $C\subseteq\p^3$.
	
	Now suppose $X$ is an ACM variety of codimension $e\geq 4$, degree $d=e+3$, and $(i)\Rightarrow (iii)$ holds for varieties of codimension $<e$. Note that $\reg X=3$ and $X$ satisfies property $\N_{2,e-2}$ by \Cref{finite_pts_Np}. Denote the image of the inner projection of $X$ from the general point $q\in X$ as $X_q\subseteq\p^{n+e-1}$.
	Then the arithmetic depth of $X_q\subseteq\p^{n+e-1}$ is same with $X\subseteq\p^{n+e}$ by \Cref{depth_inner_projection}. Hence $X_q$ is also ACM and
	\begin{align*}
		\beta_{p,1}(X) &= \beta_{p,1}(X_q)+\beta_{p-1,1}(X_q)+\dbinom{e}{p}\\
		&= p\dbinom{e}{p+1}-2\dbinom{e-1}{p-1}+(p-1)\dbinom{e}{p}-2\dbinom{e-1}{p-2}+\dbinom{e}{p}\\
		&=p\dbinom{e+1}{p+1}-2\dbinom{e}{p-1}
	\end{align*}
	for $2\leq p\leq e-2$ by \Cref{betti_ineq}. Also, we have $$\beta_{1,1}(X)=\beta_{1,1}(X_q)+\dbinom{e}{1}=\dbinom{e}{2}-2+\dbinom{e}{1}=\dbinom{e+1}{2}-2$$ and $\beta_{e-1,1}(X)=e-1$ since any $e+3$ points in general position are contained in a unique rational normal curve $C$ so that $\beta_{e-1,1}(Y)=\beta_{e-1,1}(C)=e-1$ by \Cref{finite_pts_on_rnc}.
	
	To show $(ii)\Rightarrow (i)$, we use the induction on $e$. For $e=3$, we have $\beta_{1,1}(X)=4$. Then $X$ is a variety of almost minimal degree of depth $n$ or an ACM variety of degree $d=e+3$ by \Cref{characterization_quadrics=max-2}. Hence this case is done as varieties of almost minimal degree of depth $n$ and codimension 3 have the Betti table $T(1,2)$ by \Cref{classification_Betti_VAMD}.
	Now let $X$ be a variety of codimension $e\geq 4$ with $$\beta_{p,1}(X)=p\dbinom{e+1}{p+1}-2\dbinom{e}{p-1}$$ for some integer $2\leq p\leq e-2$, and suppose the result holds for the varieties of codimension $\leq e-1$.
	Then we have
	\begin{align*}
		\beta_{p,1}(X) &\leq \beta_{p,1}(X_q)+\beta_{p-1,1}(X_q)+\dbinom{e}{p}\\
		&\leq p\dbinom{e}{p+1}-2\dbinom{e-1}{p-1}+(p-1)\dbinom{e}{p}-2\dbinom{e-1}{p-2}+\dbinom{e}{p}\\
		&= p\dbinom{e+1}{p+1}-2\dbinom{e}{p-1}
	\end{align*}
	by \Cref{main_upper_bound_not_vmd_dp}.
	Also, the equalities must hold since $$\beta_{p,1}(X)=p\dbinom{e+1}{p+1}-2\dbinom{e}{p-1}.$$ In particular, we have $$\beta_{p-1,1}(X_q)=(p-1)\dbinom{e}{p}-2\dbinom{e-1}{p-2}$$ and $$\beta_{p,1}(X_q)=p\dbinom{e}{p+1}-2\dbinom{e-1}{p-1}$$ so that $X_q$ is a variety of almost minimal degree of $\depth X=n$ and $a(X)=0$, or an ACM variety of degree $d=e+3$ by the induction hypothesis. Again by \Cref{betti_ineq}, we get
	\[
	\beta_{1,1}(X)=\beta_{1,1}(X_q)+e=\dbinom{e+1}{2}-2
	\]
	so that $X$ is a variety of almost minimal degree of depth $n$ or an ACM variety of degree $d=e+3$.
	If $X$ is a variety of almost minimal degree of depth $n$, then by Lemma \ref{T(b_1,b_2)_bound}, $X$ has the Betti table $T(1,e-1)$ so that $X$ has the Green-Lazarsfeld index $a(X)=0$.
	Hence we get $(ii)\Rightarrow (i)$.
\end{proof}

\begin{rem}
	The condition (2) in \Cref{characterizations_extremal_cases} can not be extended to
	\[
		\beta_{p,1}(X)=p\dbinom{e+1}{p+1}-2\dbinom{e}{p-1}\text{ for some integer }1\leq p\leq e-2
	\]
	by \Cref{characterization_quadrics=max-2} and Lemma \ref{T(b_1,b_2)_bound}.
	Also, it can not be extended to
	\[
		\beta_{p,1}(X)=p\dbinom{e+1}{p+1}-2\dbinom{e}{p-1}\text{ for some integer }2\leq p\leq e-1
	\]
	as there are many other examples with $\beta_{e-1,1}=e-1$.
\end{rem}

\begin{rem}\label{parallel_statement_d=e+3}
	Let $X^n\subseteq\p^{n+e}$ be a nondegenerate projective variety of codimension $e\geq 3$ and degree $d\geq e+3$.
	Then the following are equivalent:
	\begin{enumerate}
		\item $X$ is an ACM variety of degree $d=e+3$.
		\item $\beta_{p,1}(X)=p\dbinom{e+1}{p+1}-2\dbinom{e}{p-1}$ for some $1\leq p\leq e-2$.
		\item $\beta_{p,1}(X)=p\dbinom{e+1}{p+1}-2\dbinom{e}{p-1}$ for all $1\leq p\leq e-1$.
		\item $X$ satisfies $\N_{2,e-2}$ but not $\N_{2,e-1}$, and $X$ is not a complete intersection of three quadrics.\footnote{For this condition, there is an exceptional case: the complete intersections of three quadrics satisfy $\textbf{\textup{N}}_{2,1}$ but not $\textbf{\textup{N}}_{2,2}$ with the codimension 3 and degree 8. For details, \textit{c.f.} \cite{MR2862204}.}
		\item $X$ has the following Betti table
		\[\begin{array}{c|cccccc}
			&0 &1 &\cdots &e-2 &e-1 &e \\ \hline
			0&1 & & & & & \\
			1& &\beta_{1,1} &\cdots &\beta_{e-2,1} &\beta_{e-1,1} & \\
			2& & & & &e &2 
		\end{array}\]
		\noindent
		where $\beta_{p,1}=p\binom{e+1}{p+1}-2\binom{e}{p-1}$ for all $1\leq p\leq e-1$.
	\end{enumerate}
\end{rem}

\begin{rem}
	Let $X^n\subseteq\p^{n+e}$ be a nondegenerate projective variety of dimension $n\geq 1$, codimension $e$, degree $d$, and the sectional genus $\pi(X)$. Recall that the sectional genus is 0 if $X$ is a variety of minimal degree, and for the varieties of degree $d=e+2$, we have $\pi(X)\leq 1$ and
	\[
	\pi(X)=\begin{cases}
		0 &\text{if $X$ is non-ACM and $d=e+2$}\\
		1 &\text{if $X$ is ACM and $d=e+2$, i.e. $X$ is del Pezzo}.\\
	\end{cases}
	\]
	In the same way, for the varieties of degree $d=e+3$, it is easy to see $\pi(X)\leq 2$ and
	\[
	\pi(X)=\begin{cases}
		0\text{ or }1 &\text{if $X$ is non-ACM and $d=e+3$}\\
		2 &\text{if $X$ is ACM and $d=e+3$}.\\
	\end{cases}
	\]
\end{rem}

\section{Geometric descriptions of extremal varieties}

Throughout this section, let $X^n \subseteq \p^{n+e}$ be a nondegenerate projective variety of dimension $n$ and codimension $e \geq 3$ that is an extremal variety as in \Cref{characterizations_extremal_cases}. Note that $X$ is always contained in a variety $Y$ of minimal degree as a divisor. Indeed, when $d=e+2$, $X$ is contained in the rank 1 locus of a 1-generic $2\times e$ matrix by \cite[Theorem 1.4]{MR2274517}.
Also, when $d=e+3$, this follows by Green's $K_{p,1}$-theorem since $\beta_{e-1,1} (X) = e-1 >0$. Along this line, this section is devoted to describe the numerical types of $Y$ and the divisor class of $X$ in the class group of $Y$.
\\
\\
\noindent
\textit{4.1. }\textbf{Extremal varieties of almost minimal degree.}
First we consider the case $d = e+2$. The following facts are known for this case:

\begin{rem}
	Let $X$ be a variety of almost minimal degree that is not a cone. We call a rational normal scroll as \textit{an embedding scroll of $X$} if it contains $X$ as a divisor.
	\begin{itemize}
		\item When $\depth X=1$, the variety $X$ is an isomorphic projection of either a smooth rational normal scroll or the Veronese surface in $\p^5$. In the former case, $X$ is contained in a smooth rational normal scroll $Y$ as a divisor linearly equivalent to $H+2F$. 
		In the latter case, $X$ is not contained in any quadric hypersurface in $\p^4$ as its homogeneous ideal is generated by seven cubic polynomials (\cite{MR2274517}, \cite{MR2331758}).
		\item When $2\leq \depth X\leq n$, the singular locus of $X$ is a $(\depth X-2)$-dimensional linear space. Also, its embedding scroll is uniquely constructed by $J(\sing(X),X)$. Furthermore, $X$ is linearly equivalent to $H+2F$ (\cite{MR2811569}).
		\item When $\depth X=n+1$ and $X$ is smooth, the embedding scroll may not exist (\textit{e.g.} The Veronese variety $\nu_3(\p^2)$ does not lie on a threefold of minimal degree, \textit{c.f.} \cite[p. 151]{MR739785}).
		\item When $\depth X=n+1$ and $X$ is singular, the variety $Y=J(\sing(X),X)$ is an $(n+1)$-dimensional variety of minimal degree although $Y$ can be a cone over the Veronese surface in $\p^5$. In that case, $X$ must be a normal del Pezzo variety (\cite{MR1102257}, \cite{MR1162108}).
	\end{itemize}
\end{rem}

\begin{thmalphabetmaintext}\label{thm:VAMD with depth n and index 0}
	Let $X^n \subseteq \p^{n+e}$ be a nondegenerate projective variety of dimension $n$, codimension $e \geq 3$, and degree $d$. Then the following two statements are equivalent:
	\begin{enumerate}
		\item[$(i)$] X is a variety of almost minimal degree with $\depth X=n$ and $a(X)=0$.
		\item[$(ii)$] $X$ is contained in an (unique) $(n+1)$-fold rational normal scroll $Y = S(0,\ldots , 0,1,e-1)$ as a divisor linearly equivalent to $H+2F$ where $H$ and $F$ are respectively the hyperplane section and a ruling of $Y$.
	\end{enumerate}
	Furthermore, in the case above, if $X$ is not a cone, then $n \leq 4$.
\end{thmalphabetmaintext}

\begin{proof}
	$(i) \Longrightarrow (ii):$ It is enough to assume $X$ is not a cone. Since $d=e+2$ and $\depth X = n$, $\sing(X)=:\Lambda\cong \p^{n-2}$, and there is a unique rational normal $(n+1)$-fold scroll $Y = S(a_0 , a_1 ,\ldots , a_n )$ containing $X$ (\textit{c.f.} \cite[Theorem 1.3.(c)]{MR2274517}, \cite[Theorem 5.5]{MR2811569}, \cite[Theorem 1.4]{MR2274517}).
	Let $X_\Lambda$ be the image of the inner projection $$\pi_\Lambda:X\dashrightarrow \p^{e+1}$$ with the projection center $\Lambda$. Then it can be checked that $X_\Lambda$ is a smooth rational normal scroll and $Y=J(\Lambda,X_\Lambda)$.
	
	We may write $Y=S(0,\cdots,0,a,b)$ where $1\leq a\leq b$. We show that $a=1$, $b=e-1$. Let $C\subseteq\p^{e+1}$ be the general curve section of $X$. Then $C$ is a smooth curve contained in the smooth rational surface scroll $S(a,b)$. As $X$ and $C$ have the same Betti table, it follows that $a(X)=a(C)=0$. Then $a=1$ as $a(C)=a-1$ by \cite[Theorem 1.1]{MR2299577}. Hence we also have $b=e-1$. Furthurmore, $X$ is linearly equivalent to $H+2F$ (\textit{c.f.} \cite[Theorem 1.1]{MR2331758}).
	
	$(ii) \Longrightarrow (i):$ First note that $$d = H^n \cdot (H+2F) = e+2.$$ Let $C \subseteq\p^{e+1}$ be a general curve section of $X$. Then it is contained in the smooth rational normal surface scroll $S(1,e-1)$ as a divisor linearly equivalent to $H_0 +2F_0$ where $H_0$ and $F_0$ are divisors of $S(1,e-1)$. Hence $\depth C=1$ and $a(C)=0$ (\textit{c.f.} \cite[Theorem 1.1]{MR2299577}). 
	Let $\Lambda$ be the vertex of $Y$.
	As $X$ is linearly equivalent to $J(\Lambda,C)$, we have $$\depth X =\depth C+\dim \Lambda+1=n $$ (\textit{c.f.} \cite[Proposition 4.1.(2)]{MR3247023}). In particular, $a(X) =a(C)=0$.
	
	Regarding the last statement, suppose $d=e+2$, $\depth X =n$, and $X$ is not a cone. Then $X = \pi_p (\tilde{X})$ where $\tilde{X} \subseteq\p^{n+e+1}$ is a smooth rational $n$-fold scroll and
	$$\pi_p : \p^{n+e+1} \setminus \{p\} \rightarrow \p^{n+e}$$
	is the linear projection from a point $p \in \p^{n+e+1} \setminus \tilde{X}$. 	
	Let $\sec_p(\tilde{X})$ be the union of all secant lines to $X$ passing through $p$. Then the possible secant locus $\sec_p(\tilde{X})$ at $p$ is $\p^1$, $\p^2$, or $\p^3$ where $$\dim\sec_p(\tilde{X})=\depth X-1$$ (\cite[Proposition 3.2]{MR2645336}). Therefore $n=\depth X\leq 4$.
\end{proof}

We finish this subsection by providing some examples of varieties described in Theorem \ref{thm:VAMD with depth n and index 0}.

\begin{exmp}\label{example:VAMD n=1}
	For $e \geq 3$, consider the monomial curve
	$$X = \{ [  s^{e+2} : s^{e+1} t : s^{e-1} t^3 :  \cdots : st^{e+1}  : t^{e+2} ] ~|~ [s:t ]  \in \p^1 \} \subseteq\p^{e+1} .$$
	Then $X$ is obtained by an isomorphic projection of the rational normal curve $S(e+2) \subseteq\p^{e+2}$ from a point. Also the rank of the matrix
	$$\left(\begin{array}{rrrrr}
		s^{e+2}   & s^{e-1} t^3   & s^{e-2} t^{4} & \cdots & st^{e+1}  \\
		s^{e+1} t & s^{e-2} t^{4} & s^{e-3} t^{5} & \cdots & t^{e+2}
	\end{array}\right)$$
	is equal to $1$ for any $[s:t] \in \p^1$. Therefore $X$ is contained in the rational normal surface scroll
	$$Y = S(1,e-1) \subseteq\p^{e+1}$$
	which is defined as the rank $1$ locus of the matrix
	$$\left(\begin{array}{rrrrr}
		z_0 & z_2 & z_3 & \cdots & z_{e} \\
		z_1 & z_3 & z_4 & \cdots & z_{e+1}
	\end{array}\right).$$
	As a consequence, Theorem \ref{thm:VAMD with depth n and index 0} shows that $X$ is a curve of almost minimal degree with $\depth X = 1$ and $a(X) =0$.
\end{exmp}

\begin{exmp}\label{example:VAMD n=2}
	For $e \geq 3$, let
	{\small
	$$\tilde{X} = S(1,e+1)= \{ [sx:tx:s^{e+1} y : s^e t y :   \cdots : st^e y : t^{e+1} y ] ~|~ [s:t ],~[x:y] \in \p^1 \} \subseteq\p^{e+3}$$
	}
	and
	{\small
	$$X = \{ [sx:tx:s ^{e+1} y:s ^{e-1} t ^{2} y: \cdots :st ^{e} y:t ^{e+1} y]~| ~[s:t],~[x:y] \in \p^{1} \} \subseteq\p^{e+2}.$$
	}
	Since $X = \pi_p (\tilde{X})$ for $p = [0:0:0:1:0:\cdots :0] \in \p^{e+3} \setminus \tilde{X}$, we know that $X \subseteq\p^{e+2}$ is a surface of almost minimal degree. Furthermore, $\pi_p$ maps the section $S(e+1)$ of $\tilde{X}$ into a cusp curve. Thus $\pi_p \vert_{\tilde{X}} : \tilde{X} \rightarrow X$ is not an isomorphism, and hence $X$ is singular.
	
	One can observe that the rank of the matrix
	$$\left(\begin{array}{rrrrr}
		sx & s^{e-1} t^{2} y& s^{e-2} t^{3} y& \cdots & st^{e} y \\
		tx & s^{e-2} t^{3} y& s^{e-3} t^{4} y&  \cdots &t^{e+1} y
	\end{array}\right)$$
	is equal to $1$ for any $[s:t]$ and $[x:y]$ in $\p^1$. Therefore $X$ is contained in the threefold rational normal scroll
	$$Y = S(0,1,e-1) \subseteq\p^{e+2}$$
	which is defined as the rank $1$ locus of the matrix
	$$\left(\begin{array}{rrrrr}
		z_0 & z_3 & z_4 & \cdots & z_{e+1} \\
		z_1 & z_4 & z_5 & \cdots & z_{e+2}
	\end{array}\right).$$
	As a consequence, Theorem \ref{thm:VAMD with depth n and index 0} shows that $X$ is a surface of almost minimal degree with $\depth X = 2$ and $a(X) =0$.
\end{exmp}

\begin{exmp}\label{example:VAMD n=3}
	For $e \geq 3$, let $\tilde{X} = S(1,2,e-1) \subseteq\p^{e+4}$ be the smooth $3$-fold rational normal scroll parameterized as
	{\footnotesize
	$$\tilde{X} = \{ [sx:tx:s^2 y :sty : t^2 y : s^{e-1} z : s^{e-2} t z :   \cdots : t^{e-1} z ] ~|~ [s:t ] \in \p^1,~ [x:y:z]\in \p^2 \}. $$
	}
	Consider the subset
	{\footnotesize
	$$X = \{ [sx:tx:s^2 y : t^2 y : s^{e-1} z : s^{e-2} t z :   \cdots : t^{e-1} z ]~| ~[s:t ] \in \p^1,~ [x:y:z]\in \p^2 \} \subseteq\p^{e+3}.$$
	}
	Since $X = \pi_p (\tilde{X})$ for $p = [0:0:0:1:0:\cdots :0] \in \p^{e+4} \setminus \tilde{X}$, we know that $X \subseteq\p^{e+3}$ is a threefold of almost minimal degree. Furthermore, $\pi_p$ maps the section $S(2)$ of $\tilde{X}$ onto a line. Thus $\pi_p \vert _{\tilde{X}} : \tilde{X} \rightarrow X$ is not an isomorphism, and hence $X$ is singular.
	
	One can observe that $X$ is contained in the fourfold rational normal scroll
	$$Y = S(0,0,1,e-1) \subseteq\p^{e+3}$$
	which is defined as the rank $1$ locus of the matrix
	$$\left(\begin{array}{rrrrr}
		z_0 & z_4 & z_5 & \cdots & z_{e+2} \\
		z_1 & z_5 & z_6 & \cdots & z_{e+3}
	\end{array}\right).$$
	As a consequence, Theorem \ref{thm:VAMD with depth n and index 0} shows that $X$ is a threefold of almost minimal degree with $\depth X = 3$ and $a(X) =0$.
\end{exmp}

\begin{exmp}\label{example:VAMD n=4}
	For $e \geq 3$, let $\tilde{X} = S(1,1,1,e-1) \subseteq\p^{e+5}$ be the smooth $4$-fold rational normal scroll parameterized as
	{\footnotesize
	$$\tilde{X} = \{ [sx:tx:sy :ty : sz :tz : s^{e-1} w : s^{e-2} t w :   \cdots : t^{e-1} w ] ~|~ [s:t ] \in \p^1,~ [x:y:z:w]\in \p^3 \}. $$
	}
	Consider the subset
	{\footnotesize
	\begin{align*}
	    X &= \{ [sx:tx:sy -tz :ty : sz : s^{e-1} w : s^{e-2} t w :   \cdots : t^{e-1} w ]~| ~ [s:t ] \in \p^1,~ [x:y:z:w]\in \p^3 \} \\
        &\subseteq\p^{e+4}.
	\end{align*}
	}
	Since $X = \pi_p (\tilde{X})$ for $p = [0:0:1:0:0:1:0:\cdots :0] \in \p^{e+5} \setminus \tilde{X}$, we know that $X \subseteq\p^{e+4}$ is a fourfold of almost minimal degree. One can observe that $X$ is contained in the fivefold rational normal scroll
	$$Y = S(0,0,0,1,e-1) \subseteq\p^{e+4}$$
	which is defined as the rank $1$ locus of the matrix
	$$\left(\begin{array}{rrrrr}
		z_0 & z_5 & z_6 & \cdots & z_{e+3} \\
		z_1 & z_6 & z_7 & \cdots & z_{e+4}
	\end{array}\right).$$
	As a consequence, Theorem \ref{thm:VAMD with depth n and index 0} shows that $X$ is a fourfold of almost minimal degree with $\depth X = 4$ and $a(X) =0$.
\end{exmp}
~\\
\noindent
\textit{4.2. }\textbf{ACM varieties of degree=codimension+3.}
Next we consider the case $X$ is an ACM variety of degree $d=e+3$. We begin with the following lemma.

\begin{lem}\label{lem:uniqueness of embedding scroll}
	Let $X^n \subseteq \p^{n+e}$ be a nondegenerate variety of dimension $n$, codimension $e \geq 3$, and degree $d \geq e+3$. Then there exists at most one $(n+1)$-dimensional variety of minimal degree that contains $X$.
\end{lem}

\begin{proof}
	Let $Y_1 , Y_2 \subseteq \p^{n+e}$ be two $(n+1)$-dimensional varieties of minimal degree that contain $X$. Suppose that there is a point $p \in Y_1 \setminus Y_2$. We may assume that $p$ is a smooth point of $Y_1$.
	
	Let $\p^e$ be the general linear subspace of $\p^{n+e}$ passing through $p$. For each $i=1,2$, the intersection
	$$\mathcal{C}_i = Y_i \cap \p^e \subseteq \p^e$$
	is a rational normal curve of degree $e$ that contains the set $\Gamma = X \cap \p^e$ of $d$ points in linearly general position. Note that $\mathcal{C}_1 \neq \mathcal{C}_2$ since $p \in \mathcal{C}_1 \setminus \mathcal{C}_2$.
	But there exists at most one rational normal curve passing through $\Gamma$ since $d \geq e+3$. This is a contradiction. As a consequence, we can conclude that $Y_1$ and $Y_2$ must be same.
\end{proof}

\begin{remark and notation}
	Let $C$ be a projective integral curve of arithmetic genus $2$ and $\pi : \tilde{C} \rightarrow C$ be its normalization. Also let
	$$ \delta_1 (C) := \mbox{length} \left( \pi_* \mathcal{O}_{\tilde{C}} / \mathcal{O}_C \right) \quad \mbox{and} \quad \delta_2 (C) := |\mbox{Supp} \left( \pi_* \mathcal{O}_{\tilde{C}} / \mathcal{O}_C \right) |.$$
	We define \textit{the numerical type of} $C$, denoted by $\mbox{n.t}(C)$, by
	$$\mbox{n.t}(C) = \left( g(\tilde{C}),\delta_1 (C) , \delta_2 (C) \right)$$
	Since it always holds that $g(C) = g(\tilde{C})+\delta_1 (C)$, the followings are expected to occur as the numerical type of $C$:
	$$\{(2,0,0),(1,1,1),(0,2,2),(0,2,1) \}$$
	We say that a projective integral curve of arithmetic genus $g\geq 2$ is \textit{hyperelliptic} if it admits a map of degree two to the projective line. One can see that the above $C$ fails to be hyperelliptic if and only if $\mbox{n.t}(C) = (0,2,1)$.
\end{remark and notation}

\begin{cor}\label{cor:unique embedding surface of minimal degree}
	Let $C \subseteq \p^{e+1}$ be a linearly normal curve of codimension $e \geq 3$ and arithmetic genus $2$ (and hence of degree $e+3$). Then
	\begin{enumerate}[\rm (1)]
		\item If $C$ is hyperelliptic, or equivalently, $\mbox{n.t}(C) \neq (0,2,1)$, then $C$ is contained in a rational normal surface scroll $S(a,e-a)$ for some $\frac{e-3}{2} \leq a \leq \frac{e}{2}$. Furthermore, $C \equiv 2H+(3-e)F$ if $a > 0$, and $C \equiv 6R$ if $e=3$ and $a=0$ where $R \cong \p^1$ is the effective generator of the divisor class group of $S(a,e-a)$.
		\item If $\mbox{n.t}(C) = (0,2,1)$, then $C$ is contained in the singular rational normal surface scroll $S(0,e)$ and $C \equiv (e+3)R$ where $R \cong \p^1$ is the effective generator of the divisor class group of $S(0,e)$.
	\end{enumerate}
	Furthermore, the above rational normal surface scroll is the unique surface of minimal degree which contains $C$.
\end{cor}

\begin{proof}
	Since $\beta_{e-1,1} (C) = e-1 >0$, Green's $K_{p,1}$-theorem says that $C$ is contained in a surface $S$ of minimal degree as a divisor. Also when $e=4$, one can see that $S$ is not the Veronese surface since $\p^2$ contains no curve of arithmetic genus $2$. Along this line, it holds that $C$ is contained in a rational normal surface scroll $S = S(a,e-a)$ for some $0 \leq a \leq \frac{e}{2}$. The uniqueness of $S$ comes from Lemma \ref{lem:uniqueness of embedding scroll}.
	\smallskip
	
	$(1)$ Since $C$ is a hyperelliptic curve of genus $2$, there is a unique $2:1$ map $f : C \rightarrow \p^1$. Letting $A := f^* \mathcal{O}_{\p^1} (1)$, we can write
	$$\mathcal{O}_C (1) = A^m \otimes B$$
	for some effective divisor $B$ of degree $b$ on $C$ such that $h^0 (C,B )>0$ and $h^0 (C,A^{-1}\otimes B )=0$. Thus $e+3 = 2m+b$. Then it holds by \cite[Theorem 3.1.(2)]{MR2559684} that
	either $b=0,1$ and $m+b \geq 3$, or else $b = 2,3$ and $m+b \geq 4$. Also $C$ is contained in a rational normal surface scroll $S(m+b-3,m)$. Therefore $a:=m+b-3$ satisfies the desired inequalities
	$$\frac{e-3}{2} \leq a \leq \frac{e}{2}$$
	by \cite[Corollary 3.3]{MR2559684}.
	
	$(2)$ Suppose that $a \geq 1$. Then we can write $C \equiv \alpha H + \beta F$ for some $\alpha \geq 1$. If $\alpha = 1$, then $g(C) = 0$. Also if $\alpha \geq 3$, then $C$ admits a trisecant line while it is ideal-theoretically cut out by quadrics. As a consequence, $\alpha =2$, and hence $C$ is hyperelliptic which is a contradiction. Therefore, $a = 0$ and $C \equiv (e+3)R$.
\end{proof}

\begin{thmalphabetmaintext}\label{thm:ACM with degree=e+3}
	Let $X^n \subseteq \p^{n+e}$ be a nondegenerate projective variety of dimension $n$, codimension $e \geq 3$, and degree $d$. Then the following two statements are equivalent:
	\begin{enumerate}
		\item[$(i)$] $X$ is an ACM variety of degree $d=e+3$.
		\item[$(ii)$] $X$ is contained in an (unique) $(n+1)$-fold rational normal scroll $Y = S(a_0 , a_1 , \ldots ,a_n )$ for some $0 \leq a_0 \leq a_1 \leq \ldots \leq a_n$ such that either
		\begin{enumerate}
			\item[$(\alpha)$] $a_{n-1} > 0$ and the divisor class of $X$ is $2H+(3-e)F$  where $H$ and $F$ are respectively the hyperplane section and a ruling of $Y$, or
			\item[$(\beta)$] $a_{n-1} = 0$ and the divisor class of $X$ is $(e+3)R$  where $R \cong \p^n$ is the effective generator of the divisor class group of $Y$.
		\end{enumerate}
	\end{enumerate}
\end{thmalphabetmaintext}

\begin{proof}
	$(i) \Longrightarrow (ii):$ Note that $\beta_{e-1,1}(X) = e-1 > 0$ by \Cref{characterizations_extremal_cases}.
	Then Green's $K_{p,1}$-theorem says that $X$ is contained in an $(n+1)$-dimensional variety $Y$ of minimal degree.
	Let $C \subseteq\p^{e+1}$ and $S \subseteq\p^{e+1}$ be respectively a general curve section of $X$ and a general surface section of $Y$ such that $C \subseteq S$. Then $C$ is a linearly normal curve of arithmetic genus $2$, and $S$ is a surface of minimal degree. Then Corollary \ref{cor:unique embedding surface of minimal degree} shows that $C$ is contained a unique surface of minimal degree, say $S(a,e-a)$. By Corollary \ref{cor:unique embedding surface of minimal degree}, either $\frac{e-3}{2} \leq a \leq \frac{e}{2}$ or else $a=0$. This implies that $Y$ is a rational normal $(n+1)$-fold scroll, and hence it can be written as $Y = S(a_0 , a_1 , \ldots ,a_n )$ for some $0 \leq a_0 \leq a_1 \leq \ldots \leq a_n$.
	
	Now it remains to determine the divisor class of $X$ in $Y$. If $a_{n-1} \geq 1$, then the class group of $Y$ is freely generated by the hyperplane section $H$ and a ruling $F$ of $Y$. In this case, $X \equiv 2H+(3-e)F$ by Corollary \ref{cor:unique embedding surface of minimal degree}. Also, if $a_{n-1} = 0$, then $a_n = e$, and the divisor class group of $Y$ is generated by an effective divisor $R \cong \p^n$. Thus $X \equiv (e+3)R$ since $\deg X = e+3$.
	
	$(ii) \Longrightarrow (i):$ Let $C \subseteq\p^{e+1}$ and $S \subseteq\p^{e+1}$ be as above. Then $S$ is a rational normal surface scroll and $d = e+3$. When $S$ is singular, $C$ is ACM by \cite[Example 5.2]{MR1824228}. Also when $S$ is smooth, it holds that $C$ is ACM by \cite[Theorem 5.10.(a)]{MR1615938}. This implies that $d=e+3$ and $X$ is ACM.
\end{proof}

\begin{rem}\label{same_quadratic_strand}
	Let $X$ and $Y = S(0,\ldots , 0,1,e-1)$ be as in Theorem \ref{thm:VAMD with depth n and index 0}. Let $D$ be the effective divisor in the class $H+(1-e)F$. It is the join of the section $S(1)$ with the vertex of $Y$. Now, the reducible set $Z = X \cup D$ is linearly equivalent to $2H+(3-e)F$, and hence its graded Betti numbers are the same by \cite[Proposition 4.1.(1)]{MR3247023}. Indeed, $X$ is contained in $Z$ as an irreducible component, and the quadratic linear strand of the minimal free resolution of $X$ comes from that of $Z$. As a consequence, $Z$ can be interpreted as the limit of irreducible divisors in $|2H+(3-e)F|$ which are ACM and of degree $e+3$. This explains why the almost minimal degree case appear in our main theorem. It is an exceptional case comparing with the case $X$ is an ACM variety of degree $d=e+3$ in the sense that it is an irreducible component of a \textit{broken} divisor $Z$ in $|2H+(3-e)F|$.
\end{rem}

\section*{Statements and Declarations}
There are no competing interests to declare.

\bibliographystyle{amsalpha}
\bibliography{ref.bib}

\end{document}